\documentclass{article}

\usepackage{amssymb,amsmath,bm,hyperref,geometry}
\usepackage{amssymb,amsmath,bm,graphics,graphicx,url,color}
\def\no{\noindent}
\def\pmatrix{\left(\begin{array}}
\def\endpmatrix{\end{array}\right)}

\def\RR{\mathbb{R}}

\def\by{\bar{y}}

\def\I{{\cal I}}

\def\P{{\cal P}}

\def\dd{\mathrm{d}}

\newtheorem{theo}{Theorem}
\newtheorem{lem}{Lemma}
\newtheorem{cor}{Corollary}
\newtheorem{rem}{Remark}

\def\proof{\noindent\underline{Proof}\quad}
\def\QED{\mbox{~$\Box{~}$}}

\def\bfzero{{\bm{0}}}

\def\bfgamma{{\bm{\gamma}}}

\def\bfphi{{\bm{\phi}}}

\def\aa{\alpha}

\begin{document}

\title{A spectrally accurate step-by-step method for the numerical solution of fractional differential equations}

\author{
Luigi\,Brugnano\,\thanks{Dipartimento di Matematica e Informatica ``U.\,Dini'', Universit\`a di Firenze, Italy, {\tt luigi.brugnano@unifi.it}} \and
Kevin\,Burrage\,\thanks{School of Mathematical Sciences, Queensland University of Technology, Australia, {\tt kevin.burrage@qut.edu.au}} \and
Pamela\,Burrage\,\thanks{School of Mathematical Sciences, Queensland University of Technology, Australia, {\tt pamela.burrage@qut.edu.au}}\and Felice\,Iavernaro\,\thanks{Dipartimento di Matematica, Universit\`a di Bari, Italy, {\tt felice.iavernaro@uniba.it}}
}

\maketitle

\begin{abstract} In this paper we consider the numerical solution of fractional differential equations. In particular, we study a step-by-step graded mesh procedure based on an expansion of the vector field using orthonormal Jacobi polynomials. Under mild hypotheses, the proposed procedure is capable of getting spectral accuracy. A few numerical examples are reported to confirm the theoretical findings.

\bigskip
\no{\bf Keywords:} fractional differential equations, fractional integrals, Jacobi polynomials, Hamiltonian Boundary Value Methods, HBVMs.

\bigskip
\no{\bf AMS-MSC:} 65L05, 65L03, 65L99.

\end{abstract}

\section{Introduction} 

In recent years fractional differential equation modelling has become more and more frequent - see, for example, the two monographs \cite{Po1999} and \cite{Die2004}. In fact the use of fractional models has quite a long history and early applications considered their behaviour in describing anomalous diffusion in a randomly hetergeneous porous medium  \cite{LY04} and in water resources modelling \cite{BWM00}.  In mathematical biology, a number of researchers have considered fractional reaction diffusion equations - see, for example, \cite{BKB14} who have studied via a spectral approach the interfacial properties of the Allen-Cahn equation with a quartic double well potential, a model often used for alloy kinetics and the dynamics of phase boundaries. Henry, Langlands and Wearne \cite{HLW05} analysed Turing pattern formation in fractional activator-inhibitor systems.
In a computational medicine setting, Henry and Langlands \cite{HL08} developed fractional cable models for spiny neuronal dendrites. Orsingher and Behgin \cite{OB09} have considered the dynamics of fractional chemical kinetics for unimolecular systems using time change.  The authors in \cite{BKGRB14} have made a series of arguments based on Riesz potential theory and the wide range of heterogeneous cardiac tissues for the use of fractional models in cardiac electrophysiology, while \cite{CBTB15} have explored the effects of fractional diffusion on the cardiac action potential shape and electrical propagation. In addition, \cite{MFB09}, \cite{HFMTKS12} and \cite{BB17} have considered the calibration of fractional cardiac models through a fractional Bloch-Torrey equation.

Due to these application advances there is clearly a need to develop advanced numerical methods for fractional differential equations. In a purely time fractional setting the standard approach has been to use a first order approximation to the Caputo derivative. Lubich \cite{Lu1985} has extended this and developed higher order approximations based on an underlying linear multistep method that are convergent of the order of this multistep method.

There are two important issues when designing effective numerical methods for Fractional Differential equations: memory and non-locality. In the case of memory, in principle one needs to evaluate the numerical approximations all the way back to the initial time point at each time step and this becomes very computationally intensive over long time intervals. Techniques such as fixed time windowing have been proposed \cite{Po1999} and more recently a number of approaches have been considered to compute the discrete convolutions to approximate the fractional operator \cite{SLS06}, \cite{ZeTuBuKa2018}. Zeng et al. \cite{ZTB18} developed fast memory-saving time stepping methods in which the fractional operator is decoupled into a local part with fixed memory length and the history part is computed by Laguerre-Gauss quadrature, whose error is independent of the stepsize.

In the second situation the solutions to FDEs are often non-smooth and may have a strong singularity at $t = 0$. In these cases the vector field may inherit this singularity behaviour and so a constant timestep would be very inefficient. Thus graded meshes have been proposed by a number of authors \cite{LYC16,ZZK17,SOG17}. Lubich \cite{L86} on the other hand deals with singularities by introducing certain correction terms.  Other approaches have been considered such as Adams methods \cite{DiFoFr2005}, trapezoidal methods \cite{Ga2015} and Bernstein splines \cite{Sa2023} while Garrappa \cite{Ga2018} has given a survey and software tutorial on numerical methods for FDEs. 

In this paper, we consider a major improvement of the recent solution approach described in \cite{ABI2019,BI2022}, based on previous work on  Hamiltonian Boundary value Methods (HBVMs) \cite{BIT2012,BI2016,BI2018,BFCIV2022} (also used as spectral methods in time \cite{BMR2019,BMIR2019,ABI2020,ABI2023}), for solving fractional initial value problems in the form:
\begin{equation}\label{inifrac}
y^{(\aa)}(t) = f(y(t)), \qquad t\in[0,T], \qquad y(0)=y_0\in\RR^m,
\end{equation}
where, without loss of generality, we omit $t$ as a formal argument, at the r.h.s. Further, for sake of brevity, the dimension $m$ of the state vector will be also omitted, when unnecessary. 

Here, for ~$\aa\in(0,1]$, ~$y^{(\aa)}(t) \equiv D^\aa y(t)$~ is the Caputo fractional derivative:
\begin{equation}\label{Dalfa}
D^\aa g(t) = \frac{1}{\Gamma(1-\aa)} \int_0^t (t-x)^{-\aa} \left[\frac{\dd}{\dd x}g(x)\right]\dd x.
\end{equation}
The Riemann-Liouville integral associated to (\ref{Dalfa}) is given by:
\begin{equation}\label{Ialfa}
I^\aa g(t) = \frac{1}{\Gamma(\aa)}\int_0^t (t-x)^{\aa-1} g(x)\dd x.
\end{equation}
It is known that \cite{Po1999}:
\begin{equation}\label{prop}
D^\aa I^\aa g(t) = g(t), ~\quad I^\aa D^\aa g(t)=g(t)-g(0), ~\quad
I^\aa t^j = \frac{j!}{\Gamma(\aa+j+1)}t^{j+\aa}, \quad j=0,1,2,\dots.
\end{equation}
Consequently, the solution of (\ref{inifrac}) is formally given by:
\begin{equation}\label{solfrac}
y(t) = y_0 + I^\aa f(y(t)), \qquad t\in[0,T].
\end{equation}

We shall generalize the approach given in \cite{ABI2019}, by defining a step-by-step procedure, whereas the procedure given in \cite{ABI2019} was of a global nature.   This strategy, when combined with a graded mesh selection, enables better handling of singularities in the derivative of the solution to (\ref{inifrac}) at the origin, which occur, for example, when $f$ is a suitably regular function. Additionally, our approach partially mitigates issues stemming from the non-local nature of the problem.  Indeed, at a given time-step of the integration procedure, it requires the solution of a differential problem within the current time interval along with a finite sequence of pure quadrature problems that account for the contribution  brought by the history term.

With this premise, the structure of the paper is as follows: in Section~\ref{jaco} we recall some basic facts about Jacobi polynomials, for later use; in Section~\ref{qpol} we describe a {\em piecewise  quasi-polynomial approximation} to the solution of (\ref{inifrac}); in Section~\ref{error} the analysis of the new method is carried out;  in Section~\ref{num} we report a few numerical tests for assessing the theoretical achievements; at last, in Section~\ref{fine} a few concluding remarks are given.

\section{Orthonormal Jacobi polynomials}\label{jaco}
Jacobi polynomials form an orthogonal set of polynomials w.r.t. a given weighting function. In more detail, for $r,\nu>-1$:
\begin{eqnarray*}
&&\bar P^{(r,\nu)}_i(x)\in\Pi_i, 	\qquad \left(P_0^{(r,\nu)}(x)\equiv 1\right) \\
 &&\int_{-1}^1 (1-x)^r(1+x)^\nu \bar P_i^{(r,\nu)}(x)\bar P_j^{(r,\nu)}(x)\dd x ~=~ 
\frac{2^{r+\nu+1}}{2i+r+\nu+1}\frac{\Gamma(i+r+1)\Gamma(i+\nu+1)}{\Gamma(i+r+\nu+1)i!}\delta_{ij}, \\
&&\qquad i,j,=0,1,\dots,
\end{eqnarray*}
where, as is usual, $\Pi_i$ is the vector space of polynomials of degree at most $i$ and $\delta_{ij}$ is the Kronecker symbol. Consequently, the shifted and scaled Jacobi polynomials
$$P_i^{(r,\nu)}(c) ~:=~ \sqrt{(2i+r+\nu+1)\frac{\Gamma(i+r+\nu+1)i!}{\Gamma(i+r+1)\Gamma(i+\nu+1)}}\bar P_i^{(r,\nu)}(2c-1), \qquad i=0,1,\dots,$$
are orthonormal w.r.t.
$$\int_0^1(1-c)^rc^\nu P_i^{(r,\nu)}(c)P_j^{(r,\nu)}(c)\dd c ~=~ \delta_{ij}, \qquad i,j=0,1,\dots.$$
In particular, hereafter we shall consider the polynomials\,\footnote{Hereafter, for sake of brevity, we shall omit the upper indices for such polynomials.}
\begin{equation}\label{jacop}
P_i(c) ~:=~ \frac{P_i^{(\aa-1,0)}(c)}{\sqrt\aa} ~\equiv~ \frac{\sqrt{2i+\aa}}{\sqrt\aa}\bar P_i^{(\aa-1,0)}(2c-1), \qquad i=0,1,\dots,
\end{equation}
with $\aa\in(0,1]$ the same parameter in (\ref{inifrac}), such that:
\begin{equation}\label{ortoj}
P_0(c)~\equiv~ 1, \qquad \aa\int_0^1(1-c)^{\aa-1}P_i(c)P_j(c)\dd c ~=~ \delta_{ij}, \qquad i,j=0,1,\dots,
\end{equation}
where we have slightly changed the weighting function, in order that it has a unit integral:
\begin{equation}\label{unit}
\aa\int_0^1(1-c)^{\aa-1}\dd c = 1.
\end{equation}
We refer to \cite{Gau2004} and the accompanying software, for their computation, and for computing the nodes and weights of the Gauss-Jacobi quadrature of order $2k$.

\begin{rem}\label{leg}
As is clear, when in (\ref{jacop})--(\ref{unit}) $\aa=1$, we obtain the scaled and shifted Legendre polynomials, orthonormal in $[0,1]$.
\end{rem}

As is well known, the polynomials (\ref{jacop}) form an orthonormal basis for the Hilbert space $L_2[0,1]$, equipped with the scalar product
\begin{equation}\label{scalpro}
(f,g) = \aa\int_0^1 (1-c)^{\aa-1}f(c)g(c)\dd c,
\end{equation} 
and the associated norm 
\begin{equation}\label{norma}
\|f\| = \sqrt{(f,f)}.
\end{equation}
Consequently, from the Parseval identity, it follows that any square summable function can be expanded along such a basis, and the corresponding expansion is convergent to the function itself.

That said, with reference to the vector field in (\ref{inifrac}), let us consider the interval $[0,h]$, and  the expansion
\begin{equation}\label{serieh}
f(y(ch)) ~=~ \sum_{j\ge0} \gamma_j(y) P_j(c), \qquad c\in[0,1],
\end{equation}
with the Fourier coefficients given by:
\begin{equation}\label{gamj_j}
\gamma_j(y) ~=~ \aa\int_0^1 (1-c)^{\aa-1}P_j(c)f(y(ch))\dd c, \qquad j=0,1,\dots.
\end{equation}
Consequently, on the interval $[0,h]$, (\ref{inifrac}) can be rewritten as
\begin{equation}\label{inifrach}
y^{(\aa)}(ch) ~=~ \sum_{j\ge0} \gamma_j(y) P_j(c), \qquad c\in[0,1],\qquad y(0)=y_0,
\end{equation}
and, by virtue of (\ref{Ialfa})-(\ref{prop}), one obtains:
\begin{equation}\label{ych}
y(ch) ~=~ y_0+h^\aa\sum_{j\ge0} \gamma_j(y)I^\aa P_j(c), \qquad c\in[0,1].
\end{equation}
In particular, due to (\ref{ortoj}) and (\ref{unit}), one has:
\begin{equation}\label{yh}
y(h)   ~=~ y_0 + \frac{h^\aa}{\Gamma(\aa+1)}\gamma_0(y).
\end{equation}

\begin{rem}
By considering that $P_0(c)\equiv1$, from (\ref{gamj_j}) one derives that
$$\gamma_0(y) = \aa\int_0^1 (1-c)^{\aa-1} f(y(ch))\dd c.$$
Consequently, because of (\ref{Ialfa}), (\ref{yh}) becomes
$$y(h) = y_0 + I^\aa f(y(h)),$$
i.e., (\ref{solfrac}) with $t=h$. This clearly explains the use of the Jacobi basis for the expansion (\ref{serieh}).

Further, we observe that, when $\aa=1$, one retrieves the approach described in \cite{BIT2012} for ODE-IVPs.
\end{rem}

We also need the following preliminary results.

\begin{lem}\label{gamjlem}
Let \,$G:[0,h]\rightarrow V$,\, with $V$ a vector space, admit a Taylor expansion at \,$t=0$. Then, 
$$\aa\int_0^1 (1-c)^{\aa-1}P_j(c)G(ch)\dd c ~=~ O(h^j), \qquad j=0,1,\dots. $$
\end{lem}
\proof By the hypothesis and (\ref{ortoj}), one has: 
\begin{eqnarray*}\lefteqn{
\aa\int_0^1 (1-c)^{\aa-1}P_j(c)G(ch)\dd c ~=~\aa\int_0^1 (1-c)^{\aa-1}P_j(c)\sum_{\ell\ge0}\frac{G^{(\ell)}(0)}{\ell!} (ch)^\ell \dd c}\\
&=&\aa\sum_{\ell\ge0} \frac{G^{(\ell)}(0)}{\ell!} h^\ell \int_0^1 (1-c)^{\aa-1}P_j(c)c^\ell \dd c~=~ \aa\sum_{\ell\ge j} \frac{G^{(\ell)}(0)}{\ell!} h^\ell  \int_0^1 (1-c)^{\aa-1}P_j(c)c^\ell \dd c\\
&=&O(h^j).\,\QED
\end{eqnarray*}

\begin{rem} 
Concerning the above result, we observe that, by considering (\ref{scalpro})-(\ref{norma}), from the Cauchy-Schwartz theorem one derives that
$$\left|\aa \int_0^1 (1-c)^{\aa-1}P_j(c)c^\ell \dd c\right| = |(P_j,c^\ell)| \le \|P_j\| \cdot \|c^\ell\| \le 1, \qquad \forall \ell\ge j\ge0.$$
\end{rem}

\medskip
From Lemma~\ref{gamjlem}, one has the following:

\begin{cor}\label{gamj} Assume that the r.h.s. of problem (\ref{inifrac}) admits a Taylor expansion at $t=0$, then the coefficients defined at (\ref{gamj_j}) satisfy: $$\gamma_j(y) = O(h^j).$$\end{cor}

The result of the previous lemma can be weakened as follows (the proof is similar and, therefore, omitted).

\begin{lem}\label{gamjlem1}
Let \,$G:[0,h]\rightarrow V$,\, with $V$ a vector space, admit a Taylor polynomial expansion of degree $k$ with remainder  at $t=0$. Then, 
$$\aa\int_0^1 (1-c)^{\aa-1}P_j(c)G(ch)\dd c ~=~ O(h^{\min(j,k)}), \qquad j=0,1,\dots. $$
\end{lem}

In such a case, the result of Corollary~\ref{gamj} is weakened accordingly.

\begin{cor}\label{gamj-1} Assume that the r.h.s. of problem (\ref{inifrac}) admits a Taylor polynomial expansion of degree $k$ with remainder  at $t=0$. Then,  the coefficients defined at (\ref{gamj_j}) satisfy: $$\gamma_j(y) = O(h^{\min(j,k)}).$$\end{cor}

However, in general, at $t=0$ the solution may be not regular, since the derivative may be singular (this is, indeed, the case, when $f$ is a regular function). In such a case, we shall resort to the following  result.

\begin{lem}\label{halfa}
Assume that the r.h.s. in (\ref{inifrac}) is continuous for all $t\in[0,h]$. Then, for a convenient $\xi_t\in(0,t)$, one has:
$$y(t) = y_0 +  \frac{y^{(\aa)}(\xi_t)}{\Gamma(\aa+1)}t^\aa, \qquad t\in[0,h].$$
\end{lem}
\proof
In fact, from (\ref{prop}), and by using the weighted mean-value theorem for integrals, one has:
\begin{eqnarray*}
y(t)&=&y_0 + I^\aa y^{(\aa)}(t) ~=~ y_0 + \frac{1}{\Gamma(\aa)}\int_0^t (t-x)^{\aa-1} y^{(\aa)}(x)\dd x\\[2mm]
&=& y_0 + \frac{y^{(\aa)}(\xi_t)}{\Gamma(\aa)}\int_0^t (t-x)^{\aa-1}\dd x ~=~ y_0 + \frac{y^{(\aa)}(\xi_t)}{\Gamma(\aa+1)} t^\aa.\,\QED
\end{eqnarray*}

As a consequence, we derive the following weaker results concerning the Fourier coefficients (\ref{gamj_j}).

\begin{cor}\label{gamj-2}
In the hypotheses of Lemma~\ref{halfa}, and assuming $f$  is continuously differentiable in a neighborhood of $y_0$, one has:
$$\gamma_0(y) = f(y_0)+ O(h^\aa), \qquad \gamma_j(y) = O(h^\aa), \quad j=1,2,\dots.$$ 
\end{cor}
\proof
In fact, from Lemma~\ref{halfa}, one has: 
$$f(y(ch)) = f\left(y_0+\frac{y^{(\aa)}(\xi_{ch})}{\Gamma(\aa+1)}(ch)^\aa\right)=
f(y_0) + f'(y_0)\frac{y^{(\aa)}(\xi_{ch})}{\Gamma(\aa+1)}(ch)^\aa + o((ch)^\aa).$$
From this relation, and from (\ref{gamj_j}), the statement then follows.\,\QED\bigskip

For later use, we also need to discuss the quadrature error in approximating the first $s$ Fourier coefficients (\ref{gamj_j}) by means of the Gauss-Jacobi formula of order $2k$,  for a convenient $k\ge s$, whose abscissae are the zeros of $P_k$, as defined in (\ref{ortoj}), and with weights:
\begin{equation}\label{gauss} 
P_k(c_i)=0, \qquad b_i = \aa\int_0^1 (1-c)^{\aa-1}\ell_i(c)\dd c, \qquad \ell_i(c) = \prod_{j=1,\,j\ne i}^k \frac{c-c_j}{c_i-c_j}, \qquad i=1,\dots,k.
\end{equation}
That is,
\begin{equation}\label{gaussjac}
\gamma_j(y) ~\approx~ \sum_{i=1}^k b_i P_j(c_j)f(y(c_ih)) ~=:~ \hat\gamma_j(y).
\end{equation}
Concerning the quadrature errors
\begin{equation}\label{quaderr}
\Sigma_j^\aa(y,h) ~:=~ \hat\gamma_j(y)-\gamma_j(y),\qquad j=0,\dots,s-1,
\end{equation}
the following result holds true.

\begin{theo}\label{thquaderr}
Assume the function $G_{j,h}(c):=P_j(c)f(y(ch))$ be of class $C^{(2k)}[0,h]$. Then, with reference to (\ref{quaderr}), one has, for a suitable $\xi=(0,1):$
\begin{equation}\label{anal}
\Sigma_j^\aa(y,h) ~=~  K_j\frac{G_{j,h}^{(2k)}(\xi)}{(2k)!} ~=~ O(h^{2k-j}), \qquad j=0,\dots,s-1,
\end{equation}
with the constants $K_j$ independent of both $G_{j,h}$ and $h$.
\end{theo}
\proof The statement follows by considering that the quadrature is exact for polynomials in $\Pi_{2k-1}$.\,\QED\bigskip

However, as recalled above, when $f$ is a smooth function, the derivative of the solution of problem (\ref{inifrac}) may have a singularity at $t=0$. In such a case, estimates can be also derived (see, e.g., \cite{DVS1984,T2017} and \cite[Theorem\,5.1.8]{MM2008}).
However, for our purposes, because of (\ref{ortoj}), for $j=0,\dots,s-1$, from Corollary~\ref{gamj-2} we can easily derive the following one, assuming that $k\ge s$ and  (\ref{inifrach}) hold true:

\begin{eqnarray*} 
\hat\gamma_j(y) &=& \sum_{i=1}^k b_iP_j(c_i)f(y(c_ih)) ~=~ \sum_{i=1}^k b_iP_j(c_i)\sum_{\ell\ge0} P_\ell(c_i) \gamma_\ell(y)\\ \nonumber
&=& \sum_{\ell=0}^{s-1} \underbrace{\left(\sum_{i=1}^k b_i P_j(c_i)P_\ell(c_i)\right)}_{=\,\delta_{j\ell}}\gamma_\ell(y) ~+~ O(h^\aa)
~\equiv~ \gamma_j(y) + \Sigma_j(y,h), \qquad j=0,\dots,s-1.
\end{eqnarray*}
Consequently,
\begin{equation} \label{Sjyh}
\Sigma_j(y,h) ~=~ O(h^\aa), \qquad j=0,\dots,s-1.
\end{equation}

In order to derive a polynomial approximation to (\ref{inifrach}), it is enough to truncate the infinite series in (\ref{serieh}) to a finite sum, thus obtaining:

\begin{equation}\label{sums0}
\sigma^{(\aa)}(ch) ~=~ \sum_{j=0}^{s-1} \gamma_j(\sigma)P_j(c), \qquad c\in[0,1], \qquad \sigma(0)=y_0,
\end{equation}
with $\gamma_j(\sigma)$ formally given by (\ref{gamj_j}) upon replacing $y$ by $\sigma$.
In so doing, (\ref{ych}) and (\ref{yh}) respectively become:
\begin{equation}\label{sigch}
\sigma(ch) ~=~ y_0+h^\aa\sum_{j=0}^{s-1} \gamma_j(\sigma)I^\aa P_j(c), \qquad c\in[0,1],
\end{equation}
and
\begin{equation}\label{sigh}
\sigma(h)   ~=~ y_0 + \frac{h^\aa}{\Gamma(\aa+1)}\gamma_0(\sigma).
\end{equation}

It can be shown that Corollary~\ref{gamj},  Corollary~\ref{gamj-1}, Lemma~\ref{halfa}, and Corollary~\ref{gamj-2} continue formally to hold for $\sigma$. Further, by considering the approximation of the Fourier coefficients obtained by using the Gauss-Jacobi quadrature of order $2k$,
\begin{equation}\label{gaussjac-1}
\hat\gamma_j(\sigma) ~=~ \sum_{i=1}^k b_i P_j(c_j)f(\sigma(c_ih)) ~\equiv~ \gamma_j(\sigma) + \Sigma_j^\aa(\sigma,h),
\end{equation}
the result of Theorem~\ref{thquaderr} continues to hold. Moreover, we shall assume that the expansion
\begin{equation}\label{serieh1}
f(\sigma(ch)) = \sum_{j\ge0}P_j(c)\gamma_j(\sigma), \qquad c\in[0,1],
\end{equation}
holds true, similarly as (\ref{serieh}), from which also (\ref{Sjyh}) follows for the quadrature error $\Sigma_j(\sigma,h)$, when $\sigma$ has a singular derivative at $0$.   In the next sections, we shall better detail, generalize, and analyze this approach.

\section{Piecewise quasi-polynomial approximation}\label{qpol}

To begin with, in order to obtain a piecewise quasi-polynomial approximation to the solution of (\ref{inifrac}),  we consider a partition of the integration interval in the form:
\begin{equation}\label{tn}
t_n = t_{n-1} + h_n,  \qquad n=1,\dots,N, 
\end{equation}
where
\begin{equation}\label{T}
t_0=0, \qquad \sum_{n=1}^N h_n = T.
\end{equation} 
Then, according to (\ref{sums0})--(\ref{sigh}), on the first subinterval $[0,h_1]$ we can derive a polynomial approximation of degree $s-1$ to (\ref{serieh}), thus getting the fractional initial value problem
\begin{equation}\label{inifracs}
\sigma_1^{(\aa)}(ch_1) = \sum_{j=0}^{s-1} \gamma_j(\sigma_1) P_j(c), \qquad c\in[0,1], \qquad \sigma_1(0)=y_0,
\end{equation}
where $\gamma_j(\sigma_1)$ is formally given by (\ref{gamj_j}), upon replacing $y$ by $\sigma_1$ at the right-hand side: 
\begin{equation}\label{gamj1}
\gamma_j(\sigma_1) = \aa\int_0^1 (1-c)^{\aa-1}P_j(c)f(\sigma_1(ch_1))\dd c, \qquad j=0,\dots,s-1.
\end{equation}
The solution of 
(\ref{inifracs}) is a {\em quasi-polynomial} of degree $s-1+\aa$, formally given by:
\begin{equation}\label{sol}
\sigma_1(ch_1) = y_0 + h_1^\aa \sum_{j=0}^{s-1} \gamma_j(\sigma_1) I^\aa P_j(c),  \qquad c\in[0,1].
\end{equation}

However, in order to actually compute the Fourier coefficients $\{\gamma_j(\sigma_1)\}$, one needs to approximate them by using the Gauss-Jacobi quadrature (\ref{gauss}):\footnote{Hereafter, as a notational convention, $\gamma_j^i:=\hat\gamma_j(\sigma_i)$, $i=1,\dots,N$.}  
$$\gamma_j^1 ~:=~ \hat\gamma_j(\sigma_1) ~=~ \sum_{i=1}^k b_i P_j(c_i)f(\sigma_1(c_ih_1)) ~\equiv~ \gamma_j(\sigma_1) +\Sigma_j^\aa(\sigma_1,h_1).$$
In so doing,  (\ref{sol}) now formally becomes 
\begin{equation}\label{sol1}
\sigma_1(ch_1) = y_0 + h_1^\aa \sum_{j=0}^{s-1} \gamma_j^1 I^\aa P_j(c),  \qquad c\in[0,1].
\end{equation}
Moreover,  one solves the system of equations, {\em having (block) dimension $s$ independently of $k$}:
\begin{equation}\label{sist1}
\gamma_j^1 ~=~  
\sum_{i=1}^k b_i P_j(c_i)f\left( y_0+h^\aa \sum_{\nu=0}^{s-1} \gamma_\nu^1\, I^\aa P_\nu(c_i)\right), \qquad j=0,\dots,s-1.
\end{equation}
This kind of procedure is typical of HBVM$(k,s)$ methods, in the case of ODE-IVPs \cite{BI2016}: the main difference, here, derives from the non-locality of the operator. As is clear (compare with (\ref{sigh})), the approximation at $t=h_1$ will be now given by:\,\footnote{Hereafter, we shall in general denote by $\bar y_n$ the approximation to $y(t_n)$, since $y_n(t)$ will be later used to denote the restriction of $y(t)$ to the subinterval $[t_{n-1},t_n]$, $n=1,\dots,N$.}
\begin{equation}\label{sig1h}
\by_1~:=~\sigma_1(h_1) ~=~ y_0+\frac{h_1^\aa}{\Gamma(\aa+1)}\gamma_0^1.
\end{equation}

\begin{rem} For an efficient and stable evaluation of the integrals 
\begin{equation}\label{Iaci}
I^\aa P_0(c_i), ~ I^\aa P_1(c_i), ~\dots,~ I^\aa P_{s-1}(c_i), \qquad i=1,\dots,k,
\end{equation}
we refer to the procedure described in \cite{ABI2022}.
\end{rem}

\subsection{The general procedure}
We now generalize the previous procedure for the subsequent integration steps. For later use, we shall denote:
\begin{equation}\label{yn}
y_n(ch_n) ~:=~ y(t_{n-1}+ch_n), \qquad c\in[0,1],\qquad n=1,\dots,N,
\end{equation}
the restriction of the solution on the interval $[t_{n-1},t_n]$. 
Similarly, we shall denote by $\sigma(t)\approx y(t)$ the piecewise quasi-polynomial approximation such that:
\begin{equation}\label{sign}
\sigma|_{[t_{n-1},t_n]} (t_{n-1}+ch_n) ~=:~ \sigma_n(ch_n) ~\approx~ y_n(ch_n), \qquad c\in[0,1], \qquad n=1,\dots,N.
\end{equation}
Then, by using an induction argument, let us suppose one knows the quasi-polynomial approximations
\begin{equation}\label{qpoli}
\sigma_i(ch_i) = \phi^\aa_{i-1}(c,\sigma) + h_i^\aa \sum_{j=0}^{s-1} \gamma_j^i I^\aa P_j(c) ~ \approx~ y_i(ch_i), \qquad c\in[0,1],
\qquad i=1,\dots,n,
\end{equation} 
where $\phi^\aa_{i-1}(c,\sigma)$ denotes a {\em history term}, to be defined later, such that:
\begin{itemize}
\item in the first subinterval, according to (\ref{sol1}), $\sigma_0^\aa(c,\sigma)\equiv y_0$, $c\in[0,1]$;
\item for $i>1$,  $\phi_{i-1}^\aa(c,\sigma)$ only depends on  $\sigma_1,\dots,\sigma_{i-1}$. 
\end{itemize}
The corresponding  approximations at the grid-points are given by
\begin{equation}\label{yi}
\by_i~:=~\sigma_i(h_i) ~=~ \phi^\aa_{i-1}(1,\sigma)~+~\frac{h_i^\aa}{\Gamma(\aa+1)}\gamma_0^i~\approx~ y_i(h_i)~\equiv~y(t_i),\qquad i=1,\dots,n,
\end{equation}
and assume we want to compute
\begin{equation}\label{qpoln}
\sigma_{n+1}(ch_{n+1}) := \phi^\aa_n(c,\sigma)+h_{n+1}^\aa \sum_{j=0}^{s-1} \gamma_j^{n+1} I^\aa P_j(c)~\approx~y_{n+1}(ch_{n+1}), \qquad c\in[0,1].
\end{equation}
Hereafter, we shall assume that the time-steps in (\ref{tn})  define a {\em graded mesh}. In more detail, for a suitable $r>1$:
\begin{equation}\label{hn}
h_n ~=~ rh_{n-1} ~\equiv~ r^{n-1} h_1, \qquad n=1,\dots,N.
\end{equation}

\begin{rem}
As is clear, in the limit case where $r=1$, (\ref{hn}) reduces to a uniform mesh with a constant time-step\, $h=T/N$.
\end{rem}
 
 In order to derive the approximation (\ref{qpoln}), we start computing the solution of the problem in the subinterval $[t_n,t_{n+1}]$. Then, for $t\equiv t_n+ch_{n+1}$, $c\in[0,1]$, one has:
 
\begin{eqnarray*} 
\lefteqn{
y(t) ~\equiv~y_{n+1}(ch_{n+1})~=~y_0 ~+~ \frac{1}{\Gamma(\aa)}\int_0^{t_n+ch_{n+1}} (t_n+ch_{n+1}-x)^{\aa-1}f(y(x))\dd x}\\[2mm]
&=&y_0 ~+ \frac{1}{\Gamma(\aa)}\int_0^{t_n} (t_n+ch_{n+1}-x)^{\aa-1}f(y(x))\dd x~+\frac{1}{\Gamma(\aa)}\int_{t_n}^{t_n+ch_{n+1}} (t_n+ch_{n+1}-x)^{\aa-1}f(y(x))\dd x\\[2mm] 
&=&y_0 ~+~ \frac{1}{\Gamma(\aa)}\sum_{\nu=1}^n \int_{t_{\nu-1}}^{t_\nu} (t_n+ch_{n+1}-x)^{\aa-1}f(y(x))\dd x\\[2mm]  
&&~+~\frac{1}{\Gamma(\aa)}\int_{t_n}^{t_n+ch_{n+1}} (t_n+ch_{n+1}-x)^{\aa-1}f(y(x))\dd x \\[2mm]  
&=&y_0 ~+~ \frac{1}{\Gamma(\aa)}\sum_{\nu=1}^n \int_0^{h_\nu} (t_n-t_{\nu-1}+ch_{n+1}-x)^{\aa-1}f(y_\nu(x))\dd x\\[2mm]  
&&~+~\frac{1}{\Gamma(\aa)}\int_0^{ch_{n+1}} (ch_{n+1}-x)^{\aa-1}f(y_{n+1}(x))\dd x \\[2mm]  
&=& y_0 ~+~ \frac{1}{\Gamma(\aa)}\sum_{\nu=1}^n h_\nu^\aa\int_0^1 
\left(\frac{r^{n-\nu+1}-1}{r-1}+cr^{n-\nu+1}-\tau\right)^{\aa-1}f(y_\nu(\tau h_\nu))\dd \tau\\[2mm] 
&& ~+~\frac{h_{n+1}^\aa}{\Gamma(\aa)}\int_0^c (c-\tau)^{\aa-1}f(y_{n+1}(\tau h_{n+1}))\dd \tau 
~\equiv~  G_n^\aa(c,y) ~+~ I^\aa f(y_{n+1}(ch_{n+1})),
\\[2mm]  
\end{eqnarray*}
having set the {\em history term}
\begin{equation}\label{Psin}  
G_n^\aa(c,y) ~:=~y_0 ~+~ \frac{1}{\Gamma(\aa)}\sum_{\nu=1}^n h_\nu^\aa\int_0^1 
\left(\frac{r^{n-\nu+1}-1}{r-1}+cr^{n-\nu+1}-\tau\right)^{\aa-1}f(y_\nu(\tau h_\nu))\dd \tau,
\end{equation}
which is a known quantity,  since it only depends on $y_1,\dots,y_n$ (see (\ref{yn})). Consequently, we have obtained
\begin{equation}\label{yn1_0}
y_{n+1}(ch_{n+1}) =  G_n^\aa(c,y) + I^\aa f(y_{n+1}(ch_{n+1})), \qquad c\in[0,1],
\end{equation}
which reduces to (\ref{solfrac}), when $n=0$ and $t\in[0,h]$. Further, by considering the expansion
\begin{equation}\label{serienh}
f(y_{n+1}(ch_{n+1})) = \sum_{j\ge0}\gamma_j(y_{n+1}) P_j(c), \qquad c\in[0,1],
\end{equation}
with the Fourier coefficients given by:
\begin{equation}\label{gnp1}
\gamma_j(y_{n+1}) = \aa\int_0^1 (1-c)^{\aa-1}P_j(c) f(y_{n+1}(ch_{n+1}))\dd c, \qquad j\ge0,
\end{equation}
one obtains that (\ref{yn1_0}) can be rewritten as:
\begin{equation}\label{yn1_1}
y_{n+1}(ch_{n+1}) = G_n^\aa(c,y) + \frac{h_{n+1}^\aa}{\Gamma(\aa)}\sum_{j\ge 0}\gamma_j(y_{n+1})I^\aa P_j(c), \qquad c\in[0,1],
\end{equation}
with the value at $t=h$ given by:
\begin{equation}\label{yn1}
y(t_{n+1}) ~\equiv ~y_{n+1}(h_{n+1}) ~=~ G_n^\aa(1,y) + \frac{h_{n+1}^\aa}{\Gamma(\aa+1)}\gamma_0(y_{n+1}).
\end{equation}

The corresponding approximation is obtained by truncating the series in (\ref{yn1_0}) after $s$ terms, and approximating the corresponding Fourier coefficients via the Gauss-Jacobi quadrature of order $2k$:
\begin{equation}\label{sign1}
\sigma_{n+1}(ch_{n+1}) ~=~ \phi^\aa_n(c,\sigma) ~+~ \frac{h_{n+1}^\aa}{\Gamma(\aa)} \sum_{j=0}^{s-1} \gamma_j^{n+1} I^\aa P_j(c), \qquad c\in[0,1],
\end{equation}
and
\begin{equation}\label{sign1_1}
\by_{n+1} ~:=~ \sigma_{n+1}(h_{n+1}) ~=~ \phi^\aa_n(1,\sigma) ~+~ \frac{h_{n+1}^\aa}{\Gamma(\aa+1)}\gamma_0^{n+1},
\end{equation}
with $\phi_n^\aa(c,\sigma)$ an approximation of $G_n^\aa(c,y)$ in  (\ref{Psin}), defined as follows:\footnote{As anticipated above, from (\ref{finc}) it follows that $\phi_0^\aa(c,\sigma)\equiv y_0$, $c\in[0,1]$, so that (\ref{sign1}) reduces to (\ref{sol1}), whereas $\phi_n^\aa(c,\sigma)$, $n>1$, only depends on $\sigma_1,\dots,\sigma_n$.}
\begin{eqnarray}\nonumber
 \phi_n^\aa(c,\sigma)&:=&y_0 ~+~  \frac{1}{\Gamma(\aa)}\sum_{\nu=1}^n h_\nu^\aa\int_0^1 
\left(\frac{r^{n-\nu+1}-1}{r-1}+cr^{n-\nu+1}-\tau\right)^{\aa-1}\sum_{j=0}^{s-1} P_j(\tau) \gamma_j^\nu\dd \tau\\
&=&y_0 ~+~  \frac{1}{\Gamma(\aa)}\sum_{\nu=1}^n h_\nu^\aa \sum_{j=0}^{s-1}\int_0^1 \nonumber
\left(\frac{r^{n-\nu+1}-1}{r-1}+cr^{n-\nu+1}-\tau\right)^{\aa-1} P_j(\tau)\dd\tau \gamma_j^\nu\\
&\equiv&y_0 ~+~  \frac{1}{\Gamma(\aa)}\sum_{\nu=1}^n h_\nu^\aa  \sum_{j=0}^{s-1} J_j^\aa\left(\frac{r^{n-\nu+1}-1}{r-1}+cr^{n-\nu+1}\right) \gamma_j^\nu,\label{finc}
\end{eqnarray}
having set
\begin{equation}\label{Jjaell}
J_j^\aa(x) := \int_0^1(x-\tau)^{\aa-1}P_j(\tau)\dd\tau, \qquad j=0,\dots,s-1,\qquad x>1,
\end{equation}
which, as we shall see in Section~\ref{Jja}, can be accurately and efficiently computed.

As is clear, in (\ref{finc}) we have  used the approximation
\begin{equation}\label{sums1}
f(\sigma_\nu(\tau h)) ~=~ \sum_{j=0}^{s-1} P_j(\tau)\gamma_j^\nu ~+~ R_\nu^\aa(\tau),
\end{equation}
with the error term given by (see (\ref{gaussjac})-(\ref{quaderr}) and (\ref{serieh1})):
\begin{eqnarray}\nonumber
R_\nu^\aa(\tau) &:=&\sum_{j=0}^{s-1} P_j(\tau)\left[\gamma_j(\sigma_\nu)-\gamma_j^\nu\right] ~+~ \overbrace{\sum_{j\ge s} P_j(\tau)\gamma_j(\sigma_\nu)}^{=:\,E_\nu^\aa(\tau)}\\  \label{Enua}
&\equiv& \underbrace{\sum_{j=0}^{s-1} P_j(\tau) \Sigma_j^\aa(\sigma_\nu,h_\nu)}_{=:\,S_\nu^\aa(\tau)}~+~E_\nu^\aa(\tau)
~\equiv~ S_\nu^\aa(\tau) + E_\nu^\aa(\tau).
\end{eqnarray}
Because of the results of Section~\ref{jaco}, we shall hereafter assume that
\begin{equation}\label{sigEjnu}
\left| \Sigma_j^\aa(\sigma_\nu,h_\nu)\right|\le\left\{\begin{array}{cc} O(h_1^\aa), &\nu=1,\\[2mm] O(h_\nu^{2k-j}), &\nu>1,\end{array}\right.
\qquad 
\left| E_\nu^\aa(\tau)\right|\le\left\{\begin{array}{cc} O(h_1^\aa), &\nu=1,\\[2mm] O(h_\nu^s), &\nu>1.\end{array}\right.
\end{equation}
Consequently, considering that $k\ge s$, $\aa\in(0,1)$, and $j=0,\dots,s-1$, one deduces:
\begin{equation}\label{SRnu}
\left| S_\nu^\aa(\tau)\right|\le\left\{\begin{array}{cc} O(h_1^\aa), &\nu=1,\\[2mm] O(h_\nu^{s+1}), &\nu>1,\end{array}\right.
\qquad 
\left| R_\nu^\aa(\tau)\right|\le\left\{\begin{array}{cc} O(h_1^\aa), &\nu=1,\\[2mm] O(h_\nu^s), &\nu>1.\end{array}\right.
\end{equation}

In so doing, see (\ref{sign1}), the Fourier coefficients satisfy the system of equations:
\begin{equation}\label{gnp1_2}
\gamma_j^{n+1} ~=~ \sum_{\ell=1}^k b_\ell P_j(c_\ell) f_{n+1}\left(\phi_n^\aa(c_\ell,\sigma)+h_{n+1}^\aa\sum_{i=0}^{s-1} \gamma_i^{n+1} I^\aa P_i(c_\ell h)\right),
  \quad j=0,\dots,s-1,
\end{equation}
having (block)-size $s$, independently of the considered value of $k$.

Clearly,
$$\gamma_j^{n+1}~=~ \gamma_j(\sigma_{n+1}) ~+~ \Sigma_j^\aa(\sigma_{n+1},h_{n+1}),\qquad j=0,\dots,s-1$$
with $\gamma_j(\sigma_{n+1})$ formally defined as in (\ref{gnp1}), upon replacing $y_{n+1}$ with $\sigma_{n+1}$,
and $$\Sigma_j^\aa(\sigma_{n+1},h_{n+1}) = O(h_{n+1}^{2k-j}), \qquad j=0,\dots,s-1,$$ the corresponding quadrature errors.

\subsection{Solving the discrete problems}\label{solve}
The discrete problem (\ref{gnp1_2}), to be solved at each integration step, can be better cast in vector form by introducing the matrices
$$
\P_s = \pmatrix{ccc} P_0(c_1) & \dots & P_{s-1}(c_1)\\ \vdots & &\vdots\\ P_0(c_k) & \dots & P_{s-1}(c_k)\endpmatrix,
\quad \I_s^\aa = \pmatrix{ccc} I^\aa P_0(c_1)& \dots &I^\aa P_{s-1}(c_1)\\ \vdots & &\vdots\\ I^\aa P_c(c_k) & \dots &I^\aa P_{s-1}(c_k)\endpmatrix\quad \in\RR^{k\times s},$$
$$\Omega = \pmatrix{ccc} b_1\\ &\ddots\\ &&b_k\endpmatrix,$$
and the (block) vectors:
$$\bfgamma^{n+1} = \pmatrix{c} \gamma_0^{n+1}\\ \vdots\\ \gamma_{s-1}^{n+1}\endpmatrix\in\RR^{sm},\qquad 
\bfphi_n^\aa = \pmatrix{c} \phi_n^\aa(c_1,\sigma)\\ \vdots \\ \phi_n^\aa(c_k,\sigma)\endpmatrix\in\RR^{km}.$$
In fact, in so doing (\ref{gnp1_2}) can be rewritten as:\footnote{As is usual, the function $f$, here evaluated in a (block) vector of dimension $k$, denotes the (block) vector made up by $f$ evaluated in each (block) entry of the input argument.}
\begin{equation}\label{vform}
\bfgamma^{n+1} = \P_s^\top\Omega \otimes I_m f\left( \bfphi_n^\aa +h_{n+1}^\aa\I_s^\aa\otimes I_m\bfgamma^{n+1}\right),
\end{equation}
This formulation is very similar to that used for HBVMs in the case $\aa=1$ \cite{BIT2011}. The formulation (\ref{vform}) has a twofold use:
\begin{itemize}
\item it shows that, assuming, for example, $f$ Lipschitz continuous, then the solution exists and is unique, for all sufficiently small timesteps $h_{n+1}$;

\item it induces a straightforward fixed-point iteration:
\begin{eqnarray}\label{itera}
\bfgamma^{n+1,0}&=&\bfzero,\\[1mm] \nonumber
\bfgamma^{n+1,\ell} &=& \P_s^\top\Omega \otimes I_m f\left( \bfphi_n^\aa +h_{n+1}^\aa\I_s^\aa\otimes I_m\bfgamma^{n+1,\ell-1}\right),\qquad \ell =1,2,\dots,
\end{eqnarray}
which we shall use for the numerical tests.
\end{itemize}

\no In fact, the following result holds true.

\begin{theo} Assume $f$ be Lipchitz with constant $L$ in in the interval $[t_n,t_{n+1}]$. Then, the iteration (\ref{itera}) is convergent for all timesteps $h_{n+1}$ such that 
\begin{equation}\label{contraz}
h_{n+1}^\aa L\|\P_s^\top\Omega\|\|\I_s^\aa\|<1.
\end{equation}
\end{theo}

\proof In fact, one has:
\begin{eqnarray*}
\lefteqn{\|\bfgamma^{n+1,\ell+1}-\bfgamma^{n+1,\ell}\|}\\
&=& \|\P_s^\top\Omega \otimes I_m \left[ f\left( \bfphi_n^\aa +h_{n+1}^\aa\I_s^\aa\otimes I_m\bfgamma^{n+1,\ell}\right)-f\left( \bfphi_n^\aa +h_{n+1}^\aa\I_s^\aa\otimes I_m\bfgamma^{n+1,\ell-1}\right)\right]\|\\
&\le& h_{n+1}^\aa L\|\P_s^\top\Omega\|\|\I_s^\aa\| \cdot \|\bfgamma^{n+1,\ell}-\bfgamma^{n+1,\ell-1}\|,
\end{eqnarray*}
hence the iteration function defined at (\ref{itera}) is a contraction, when (\ref{contraz}) holds true.\,\QED\bigskip

A simplified Newton-type iteration, akin to that defined in \cite{BIT2011} for HBVMs, will be the subject of future investigations.

\begin{rem}
We observe that the discrete problem (\ref{vform}) can be cast in a Runge-Kutta type form. In fact, the vector
\begin{equation}\label{Yn1}
Y^{n+1} ~:=~ \bfphi_n^\aa +h_{n+1}^\aa\I_s^\aa\otimes I_m\bfgamma^{n+1},
\end{equation}
in argument to the function $f$ at the r.h.s. in (\ref{vform}), can be regarded as the stage vector of a Runge-Kutta method, tailored for the problem at hand. Substituting (\ref{vform}) into (\ref{Yn1}), and using $Y^{n+1}$ as argument of $f$, then gives
\begin{equation}\label{stage}
Y^{n+1} ~=~ \phi_n^\aa +h_{n+1}^\aa\I_s^\aa\P_s^\top\Omega \otimes I_m f(Y^{n+1}).
\end{equation}
It is worth noticing that (\ref{stage}) has (block) dimension $k$, instead of $s$. However, by considering that usually $k\gg s$ (see Section~\ref{num}), it follows that solving (\ref{stage}) is less efficient than solving (\ref{vform}). 

This fact is akin to what happens for a HBVM$(k,s)$ method \cite{BIT2012,BI2016,BI2018}, which is the $k$-stage Runge-Kutta method obtained when $\aa=1$. In fact, for such a method, the discrete problem can also be cast in the form (\ref{vform}), then having (block) dimension $s$, independently of $k$ (which is usually much larger than $s$).
\end{rem}

\subsection{Computing the integrals $J_j^\aa(x)$}\label{Jja}
\begin{figure}
\centering
\includegraphics[width=9cm]{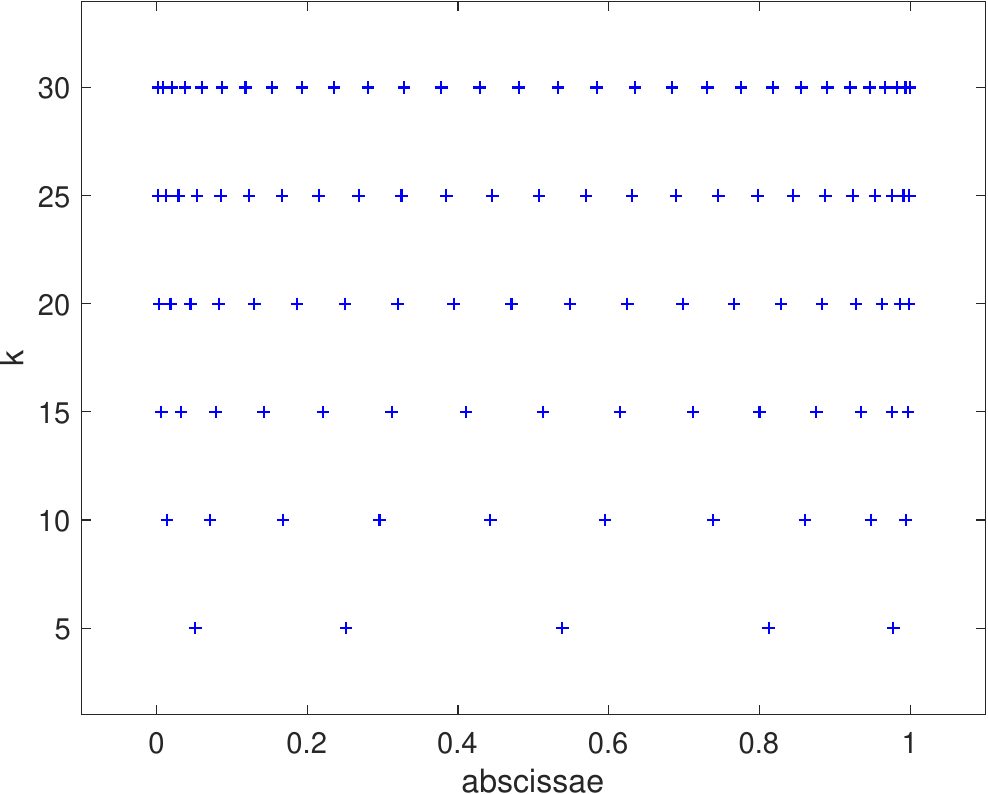}
\caption{Gauss-Jacobi abscissae for $\aa=0.5$.}
\label{absci}
\end{figure}

From (\ref{gnp1_2}) and (\ref{finc}), it follows that one needs to compute the integrals
\begin{equation}\label{tutti}
J_j^\aa\left(\frac{r^{n-\nu+1}-1}{r-1}+c_ir^{n-\nu+1}\right), \qquad i=1,\dots,k,\qquad \nu=1,\dots,n,
\end{equation}
with $\{c_1,\dots,c_k\}$ the abscissae of the Gauss-Jacobi quadrature (\ref{gauss}).  As an example, in Figure~\ref{absci} we plot the Gauss-Jacobi abscissae for $\aa=0.5$ and $k=5,10,15,20,25,30$.

Further, there is numerical evidence that a sufficiently high-order Gauss-Legendre formula can compute the integrals (\ref{Jjaell})  up to round-off, when $x\ge1.5$. Namely,   
\begin{equation}\label{nultn}
J_j^\aa\left(\frac{r^{n-\nu+1}-1}{r-1}+c_ir^{n-\nu+1}\right), \qquad i=1,\dots,k,\qquad \nu=1,\dots,n-1,
\end{equation}
and
\begin{equation}\label{nueqn0}
J_j^\aa\left(1+c_ir\right), \qquad{\rm s.t.}\qquad c_ir\ge0.5, \qquad  j=0,\dots,s-1. 
\end{equation}

However, this is no more the case, when computing:
\begin{equation}\label{nueqn}
J_j^\aa\left(1+c_ir\right), \qquad{\rm s.t.}\qquad c_ir<0.5, \qquad  j=0,\dots,s-1.
\end{equation}
Consequently, there is the problem of accurately and efficiently computing these latter integrals.

\begin{rem}
We observe that the evaluation of the integrals (\ref{nultn})-(\ref{nueqn0}),  $j=0,\dots,s-1$, via a $2p$-order Gauss-Legendre formula is inexpensive. As matter of fact, only the values $P_j(\xi_i)$ need to be computed, $j=0,\dots,s-1$, with
\begin{equation}\label{GLp}
\{\xi_1, \dots, \xi_p \}
\end{equation}
the Gauss-Legendre abscissae of the considered formula. Further, we observe that, with reference to (\ref{tutti}), only the integrals
$$J_j^\aa\left(\frac{r^n-1}{r-1}+c_ir^{n}\right),\qquad j=0,\dots,s-1, \qquad i=1,\dots,k,$$
need to be actually computed: as matter of fact, the remaining integrals,
$$J_j^\aa\left(\frac{r^\nu-1}{r-1}+c_ir^\nu\right),\qquad j=0,\dots,s-1, \qquad i=1,\dots,k,\qquad \nu=1,\dots,n-1,$$
are inherited from the previous steps.
\end{rem}

The starting point to derive an efficient algorithm for computing the integrals (\ref{nueqn}), is that the   Jacobi polynomials (\ref{jacop}) satisfy, as does any other family of orthogonal polynomials, a three-term recurrence:
\begin{equation}\label{3term}
P_{j+1}(c) = (a_jc-b_j)P_j(c) - d_jP_{j-1}(c), \quad j=0,1,\dots, \qquad P_0(c)\equiv1, \quad P_{-1}(c)\equiv 0,
\end{equation}
with prescribed $a_j,b_j,d_j$, $j\ge0$. In fact,  by setting  $\phi(c)=c$, and using the scalar product (\ref{scalpro}), to enforce (\ref{ortoj}), for $P_0,\dots,P_{s-1}$, one obtains:
\begin{equation}\label{abd}
a_j = \frac{1}{(P_{j+1},\phi\cdot P_j)}, \qquad b_j = \frac{(P_j,\phi\cdot P_j)}{(P_{j+1},\phi\cdot P_j)},\qquad d_j = \frac{(P_j,\phi\cdot P_{j-1})}{(P_{j+1},\phi\cdot P_j)}, \qquad j=0,\dots,s-2.
\end{equation}
Consequently, by recalling the definition (\ref{Jjaell}), from (\ref{3term}) one has:
$$J_0^\aa(x) = \frac{x^\aa - (x-1)^\aa}{\aa}, \qquad J_{-1}^\aa(x) = 0, \qquad \aa>0,\qquad x>1,$$ and
\begin{eqnarray*}
J_{j+1}^\aa(x) &=& \int_0^1 (x-\tau)^{\aa-1} P_{j+1}(\tau)\dd\tau\\
&=& a_j \int_0^1 (x-\tau)^{\aa-1}\tau P_j(\tau)\dd\tau - b_j J_j^\aa(x) - d_j J_{j-1}^\aa(x)\\
&=& -a_j\int_0^1(x-\tau)^\aa P_j(\tau)\dd\tau +[a_jx-b_j]J_j^\aa(x) - d_j J_{j-1}^\aa(x)\\[2mm]
&=& -a_jJ_j^{\aa+1}(x)+[a_jx-b_j]J_j^\aa(x) - d_j J_{j-1}^\aa(x), \qquad j=0,\dots,j-2.
\end{eqnarray*}

From the last two formulae, one derives that {\tt Jjalfa(a,b,d,alfa,1+ci*r)} computes all the integrals in (\ref{nueqn}) corresponding to the abscissa {\tt ci}, with {\tt Jjalfa} the Matlab$^\copyright$ function listed in Table~\ref{algo}, and the vectors {\tt a,b,d} containing the scalars (\ref{abd}).  An implementation of the previous function, using the standard variable precision arithmetic (i.e., using {\tt vpa} in Matlab$^\copyright$), allows handling values of $s$ up to 20, at least.

\begin{table}[t]
\caption{\label{algo} Matlab$^\copyright$ function {\tt Jjalfa}}
\hrulefill
\begin{verbatim}
function Ialfa = Jjalfa( a, b, d, alfa, x )
%
% Matlab function for computing the integrals J_j^alfa(x), j=0...s-1.
%
s1     = length(a);  % s1 == s-1
Ialfa  = zeros(s1+1,1);
valfa  = alfa+(0:s1);
I1     = ( x.^valfa -(x-1).^valfa )./valfa;
Ialfa(1) = I1(1);
if s1>=1
   I2 = ( a(1)*x -b(1) )*I1(1:s1) -a(1)*I1(2:s1+1);
   Ialfa(2) = I2(1); 
   for j = 2:s1
       I0 = I1; I1 = I2;
       I2 = ( a(j)*x -b(j) )*I1(1:s1-j+1) -a(j)*I1(2:s1-j+2) -d(j)*I0(1:s1-j+1);
       Ialfa(j+1) = I2(1); 
   end
end   
return
\end{verbatim}
\hrulefill
\end{table}

\section{Analysis of the method}\label{error}
From (\ref{yn1_1}) and (\ref{sign1}), and with reference to (\ref{sums1})--(\ref{SRnu}), one derives:

\begin{eqnarray*}
\lefteqn{y_{n+1}(ch_{n+1})-\sigma_{n+1}(ch_{n+1})}\\
&=& \left[G_n^\aa(c,y)-\phi_n^\aa(c,\sigma)\right] ~+~ \frac{h_{n+1}^\aa}{\Gamma(\aa)}\sum_{j=0}^{s-1}
\left[ \gamma_j(y_{n+1})-\gamma_j^{n+1}\right] I^\aa P_j(c)\\
&& ~+~  \frac{h_{n+1}^\aa}{\Gamma(\aa)}\sum_{j\ge s}\gamma_j(y_{n+1}) I^\aa P_j(c)\\
&=& \frac{1}{\Gamma(\aa)}\sum_{\nu=1}^n h_\nu^\aa\int_0^1 
\left(\frac{r^{n-\nu+1}-1}{r-1}+cr^{n-\nu+1}-\tau\right)^{\aa-1}\left[ f(y_\nu(\tau h_\nu))- f(\sigma_\nu(\tau h_\nu) +R_\nu^\aa(\tau)\right]\dd \tau\\
&&+\frac{h_{n+1}^\aa}{\Gamma(\aa)}\sum_{j=0}^{s-1}
\left[ \gamma_j(y_{n+1})-\gamma_j(\sigma_{n+1}) +\Sigma_j(\sigma_{n+1},h_{n+1})\right] I^\aa P_j(c)~+~
\frac{h_{n+1}^\aa}{\Gamma(\aa)}E_{n+1}^\aa(c).
\end{eqnarray*}

Assuming $f$ Lipschitz with constant $L$  in a suitable neighborhood of the solution, and setting
\begin{eqnarray*}
R_\nu^\aa &:=& \max_{\tau\in[0,1]} |R_\nu^\aa(\tau)|,\qquad \nu=1,\dots,N-1,\\
E_{n+1}^\aa &:=&  \max_{c\in[0,1]} |E_{n+1}^\aa(c)|,\qquad n=1,\dots,N-1,\\
e_n &:=& \max_{c\in[0,1]}|y_n(ch_n)-\sigma_n(ch_n)|, \qquad n=1,2,\dots,N,
\end{eqnarray*}
one then obtains that, for $h_{n+1}^\aa$ sufficiently small, and a suitable constant $K_1>0$:
\begin{eqnarray*}
e_{n+1} &\le&\frac{K_1L}{\Gamma(\aa)}\sum_{\nu=1}^n h_\nu^\aa\int_0^1 
\left(\frac{r^{n-\nu+1}-1}{r-1}-\tau\right)^{\aa-1}\dd\tau e_\nu\\[1mm]
&&~+~\frac{K_1}{\Gamma(\aa)}\sum_{\nu=1}^n h_\nu^\aa\int_0^1 
\left(\frac{r^{n-\nu+1}-1}{r-1}-\tau\right)^{\aa-1}\dd\tau R_\nu^\aa 
~+~g_{n+1},\qquad n=1,\dots,N-1,
\end{eqnarray*}
with 
\begin{equation}\label{gn1a}
g_{n+1}^\aa = O(h_{n+1}^{s+\aa}). 
\end{equation}

Considering that
\begin{eqnarray}\nonumber
\lefteqn{\frac{1}{\Gamma(\aa)}\int_0^1 \left(\frac{r^{n-\nu+1}-1}{r-1}-\tau\right)^{\aa-1}\dd\tau} \\ \label{sum2}
&=&\frac{1}{\Gamma(\aa+1)}\left[ \left(\frac{r^{n-\nu+1}-1}{r-1}\right)^\aa-\left(\frac{r^{n-\nu+1}-r}{r-1}\right)^\aa\right]~\le~\frac{1}{\Gamma(\aa+1)},
\qquad \nu=1,\dots,n,~
\end{eqnarray}
and (see (\ref{SRnu}))
$$h_1^\aa R_1^\aa \le O(h_1^{2\aa}), \qquad h_\nu^\aa R_\nu^\aa \le O(h_\nu^{s+\aa}), \quad \nu>1,$$
so that (see (\ref{hn})), for a suitable constant $K_2>0$,
\begin{eqnarray*}
\frac{K_1}{\Gamma(\aa+1)}\sum_{\nu=1}^n h_\nu^\aa R_\nu^\aa &\le& K_2\left( h_1^{2\aa} + \sum_{\nu=2}^n h_\nu^{s+\aa}\right)
~=~K_2\left( h_1^{2\aa} + h_n^{s+\aa}\sum_{\nu=0}^{n-2} (r^{s+\aa})^{-\nu}\right)\\
&\le& K_2\left( h_1^{2\aa} + \frac{h_n^{s+\aa}}{1-r^{-(s+\aa)}}\right),
\end{eqnarray*}
one eventually obtains:
$$
e_{n+1} ~\le~\frac{K_1L}{\Gamma(\aa)}\sum_{\nu=1}^n h_\nu^\aa\int_0^1 
\left(\frac{r^{n-\nu+1}-1}{r-1}-\tau\right)^{\aa-1}\dd\tau\, e_\nu~+~\psi_n^\aa, \qquad n=1,\dots,N-1,
$$
with (see (\ref{hn}) and (\ref{gn1a}))
$$\psi_n^\aa ~:=~ K_2\left( h_1^{2\aa} + \frac{h_n^{s+\aa}}{1-r^{-(s+\aa)}}\right) + g_{n+1}^\aa ~=~ O(h_1^{2\aa} + h_{n+1}^{s+\aa}), \qquad n=1,\dots,N-1.$$
At last, by considering that, for $n=1,\dots,N-1$,
\begin{eqnarray}\nonumber
\lefteqn{\frac{1}{\Gamma(\aa)}\sum_{\nu=1}^n h_\nu^\aa\int_0^1 
\left(\frac{r^{n-\nu+1}-1}{r-1}-\tau\right)^{\aa-1}\dd\tau}\\ \nonumber
&=&\frac{1}{\Gamma(\aa+1)} \sum_{\nu=1}^n\left[ h_\nu^\aa 
\left(\frac{r^{n-\nu+1}-1}{r-1}\right)^\aa-h_{\nu+1}^\aa 
\left(\frac{r^{n-\nu}-1}{r-1}\right)^\aa\right]\\ \label{sum1}
&=&\frac{1}{\Gamma(\aa+1)}\left[ \left(h_1\frac{r^n-1}{r-1}\right)^\aa -h_{n+1}^\aa\right]
~\le~\frac{T^\aa}{\Gamma(\aa+1)},
\end{eqnarray}
setting 
$$\rho = \exp\left(\frac{K_1LT^\aa}{\Gamma(\aa+1)}\right),$$
and considering that $e_1\le O(h_1^{2\aa})$, from the Gronwall lemma  (see, e.g.,  \cite{LaTri1988}) one eventually obtains, for a suitable $K>0$:
\begin{eqnarray}\nonumber
e_{n+1} &\le&K\rho\left( h_1^{2\aa} + \sum_{\nu=1}^n \psi_n^\aa \right) ~=~ K\rho\left( (n+1)h_1^{2\aa} + \frac{h_{n+1}^{s+\aa}}{1-r^{-(s+\aa)}}\right)\\ \label{bound}
&\le& K\rho\left( Nh_1^{2\aa} + \frac{h_N^{s+\aa}}{1-r^{-(s+\aa)}}\right), \qquad n=1,\dots,N-1.
\end{eqnarray}
Considering that (see (\ref{T}) and (\ref{hn})) $N = \log_r\left( 1 + T(r-1)/h_1\right)$, the error is optimal when $h_1^{2\aa}\sim h_N^{s+\aa}$.

\begin{rem} In the case where the vector field of problem (\ref{inifrac}) is suitably smooth at $t=0$, so that a constant timestep $h=T/N$ can be used, the estimate
\begin{equation}\label{heq}
e_n \le O(nh^{s+\aa}), \qquad n=1,\dots,N,
\end{equation}
can be derived for the error, by using similar arguments as above.
\end{rem}

\section{Numerical Tests}\label{num}
We here report a few numerical tests to illustrate the theoretical findings. For all tests, when not otherwise specified, we use $k=30$ for the Gauss-Jacobi quadrature (\ref{gauss}), and $p=30$ for the Gauss-Legendre quadrature formula (\ref{GLp}). Also, the following values of $s$ will be considered: $$s\in\{1,2,3,4,5,6,7,8,9,10,20\}.$$ All the numerical tests have been performed in Matlab$^\copyright$.  As anticipated, we use the fixed-point iteration (\ref{itera}) to solve the generated discrete problems (\ref{vform}). The iteration is carried out until full machine accuracy is gained.

The first problem, taken from \cite{Ga2018}, is given by:

\begin{equation}\label{prob1}
y^{(0.6)} = -10 y, \qquad t\in[0,5], \qquad y(0) = 1,
\end{equation}
whose solution is given by the following Mittag-Leffler function:\footnote{We refer to \cite{Ga2015ml} and the accompanying software, for its efficient Matlab$^\copyright$ implementation.}
$$E_{0.6}(t) = \sum_{j\ge0} \frac{(-10 t^{0.6})^j}{\Gamma(0.6j+1)}.$$
In Table~\ref{tab1}, we list the obtained results, in terms of maximum error, when using different values of $s$ and $h_1$. We use the value $r=1.01$ in (\ref{hn}), and the number of timesteps $N$ is chosen in order that $t_N$ is the closest point to the end of the integration interval. 
The {\tt ***}, in the line corresponding to $s=1$, means that the solution is not properly evaluated: this is due to the failure of the fixed-point iteration (\ref{itera}) in the last integration steps. As is expected from (\ref{bound}), the error decreases as $s$ increases and the initial timestep $h_1$ decreases.
Further, having a large value of $s$ is  not effective, if $h_1$ is not suitably small, and vice versa (again from the bound (\ref{bound})). 
Remarkably, by choosing $s$ large enough and $h_1$ suitably small, full machine accuracy is  gained (cf. the last 4 entries in the last column of the table, having the same error). 

In Table~\ref{tab2} we list the results obtained by using $k=s$, which is the minimum value allowed. In such a case, the accuracy is generally slightly worse, due to the fact that the quadrature error is of the same order as the truncation error. It is, however, enough choosing $k$ only slightly larger than $s$, to achieve a comparable accuracy, as is shown in Table~\ref{tab3} for $k=s+5$. Nevertheless, it must be stressed that choosing larger values of $k$ is not an issue, since the discrete problem (\ref{vform}), to be solved at each integration step, has (block) size $s$, independently of $k$.

\begin{table}
\caption{\label{tab1} Maximum error for Problem (\ref{prob1}), $r=1.01$ and $k=30$.}

\smallskip
\centering

{\tt
\begin{tabular}{|r|r r r r r r|}
\hline
 $h_1$ & $10^{-4}$ & $10^{-5}$ & $10^{-6}$ & $10^{-7}$ & $10^{-8}$ & $10^{-9}$ \\ 
 \hline
 $s\,\backslash\,N$    & 626     &    857   &     1088  &      1320  &      1551  &      1783\\
 \hline
 1~~~~ &           *** &           *** &           *** &           *** &           *** &           ***  \\ 
 2~~~~ & 3.73e-06 & 2.43e-07 & 5.40e-08 & 5.38e-08 & 5.37e-08 & 5.37e-08 \\ 
 3~~~~ & 5.81e-07 & 3.78e-08 & 2.40e-09 & 1.52e-10 & 4.64e-11 & 4.64e-11 \\ 
 4~~~~ & 1.47e-07 & 9.60e-09 & 6.11e-10 & 3.86e-11 & 2.44e-12 & 1.56e-13 \\ 
 5~~~~ & 5.29e-08 & 3.46e-09 & 2.20e-10 & 1.39e-11 & 8.79e-13 & 5.37e-14 \\ 
 6~~~~ & 2.04e-08 & 1.33e-09 & 8.49e-11 & 5.37e-12 & 3.41e-13 & 2.32e-14 \\  
 7~~~~ & 1.19e-08 & 7.81e-10 & 4.98e-11 & 3.15e-12 & 1.97e-13 & 1.08e-14 \\ 
 8~~~~ & 4.26e-09 & 2.77e-10 & 1.76e-11 & 1.11e-12 & 7.18e-14 & 7.91e-15 \\ 
 9~~~~ & 4.56e-09 & 3.03e-10 & 1.94e-11 & 1.22e-12 & 7.57e-14 & 7.91e-15 \\ 
10~~~~ & 1.89e-09 & 1.22e-10 & 7.79e-12 & 4.90e-13 & 2.96e-14 & 7.91e-15 \\ 
20~~~~ & 1.84e-09 & 1.22e-10 & 7.77e-12 & 4.90e-13 & 2.96e-14 & 7.91e-15 \\ 
\hline
\end{tabular}}
\end{table}
\begin{table}
\caption{\label{tab2} Maximum error for Problem (\ref{prob1}), $r=1.01$ and $k=s$.}

\smallskip
\centering

{\tt
\begin{tabular}{|r|r r r r r r|}
\hline
 $h_1$ & $10^{-4}$ & $10^{-5}$ & $10^{-6}$ & $10^{-7}$ & $10^{-8}$ & $10^{-9}$ \\ 
 \hline
 $s\,\backslash\,N$    & 626     &    857   &     1088  &      1320  &      1551  &      1783\\
 \hline
 1~~~~ &           *** &           *** &           *** &           *** &           *** &           ***  \\ 
 2~~~~ & 7.52e-06 & 4.99e-07 & 5.78e-08 & 5.76e-08 & 5.76e-08 & 5.76e-08 \\
 3~~~~ & 2.27e-06 & 1.51e-07 & 9.68e-09 & 6.13e-10 & 4.55e-11 & 4.55e-11 \\ 
 4~~~~ & 9.67e-07 & 6.41e-08 & 4.09e-09 & 2.59e-10 & 1.64e-11 & 1.03e-12 \\ 
 5~~~~ & 4.91e-07 & 3.26e-08 & 2.08e-09 & 1.32e-10 & 8.33e-12 & 5.23e-13 \\ 
 6~~~~ & 2.82e-07 & 1.87e-08 & 1.19e-09 & 7.56e-11 & 4.77e-12 & 2.99e-13 \\  
 7~~~~ & 1.76e-07 & 1.16e-08 & 7.44e-10 & 4.71e-11 & 2.97e-12 & 1.85e-13 \\ 
 8~~~~ & 1.16e-07 & 7.71e-09 & 4.92e-10 & 3.12e-11 & 1.97e-12 & 1.22e-13 \\ 
 9~~~~ & 8.07e-08 & 5.35e-09 & 3.42e-10 & 2.16e-11 & 1.36e-12 & 8.42e-14 \\ 
10~~~~ & 5.82e-08 & 3.86e-09 & 2.46e-10 & 1.56e-11 & 9.83e-13 & 6.01e-14 \\
20~~~~ & 6.63e-09 & 4.39e-10 & 2.80e-11 & 1.77e-12 & 1.11e-13 & 7.80e-15 \\ 
\hline
\end{tabular}}
\end{table}
\begin{table}
\caption{\label{tab3} Maximum error for Problem (\ref{prob1}), $r=1.01$ and $k=s+5$.}

\smallskip
\centering

{\tt
\begin{tabular}{|r|r r r r r r|}
\hline
 $h_1$ & $10^{-4}$ & $10^{-5}$ & $10^{-6}$ & $10^{-7}$ & $10^{-8}$ & $10^{-9}$ \\ 
 \hline
 $s\,\backslash\,N$    & 626     &    857   &     1088  &      1320  &      1551  &      1783\\
 \hline
 1~~~~ &           *** &           *** &           *** &           *** &           *** &           ***  \\ 
 2~~~~ & 3.58e-06 & 2.33e-07 & 5.40e-08 & 5.38e-08 & 5.38e-08 & 5.37e-08 \\ 
  3~~~~ &  6.56e-07 & 4.29e-08 & 2.73e-09 & 1.73e-10 & 4.63e-11 & 4.63e-11 \\
 4~~~~ &  8.57e-08 & 5.43e-09 & 3.43e-10 & 2.16e-11 & 1.36e-12 & 8.76e-14 \\
 5~~~~ &   9.20e-08 & 6.12e-09 & 3.92e-10 & 2.48e-11 & 1.56e-12 & 9.70e-14 \\
 6~~~~ &   4.38e-08 & 2.87e-09 & 1.83e-10 & 1.16e-11 & 7.30e-13 & 4.41e-14 \\ 
 7~~~~ &   3.39e-08 & 2.28e-09 & 1.46e-10 & 9.25e-12 & 5.83e-13 & 3.51e-14 \\
 8~~~~ &   2.58e-08 & 1.70e-09 & 1.09e-10 & 6.87e-12 & 4.32e-13 & 2.53e-14 \\
 9~~~~ &   2.03e-08 & 1.35e-09 & 8.60e-11 & 5.44e-12 & 3.43e-13 & 1.97e-14 \\
10~~~~ &  1.64e-08 & 1.09e-09 & 6.93e-11 & 4.38e-12 & 2.76e-13 & 1.54e-14 \\
20~~~~ &   3.28e-09 & 2.17e-10 & 1.38e-11 & 8.74e-13 & 5.40e-14 & 8.19e-15 \\
\hline
\end{tabular}}
\end{table}

\begin{figure}
\centering
\includegraphics[width=12cm]{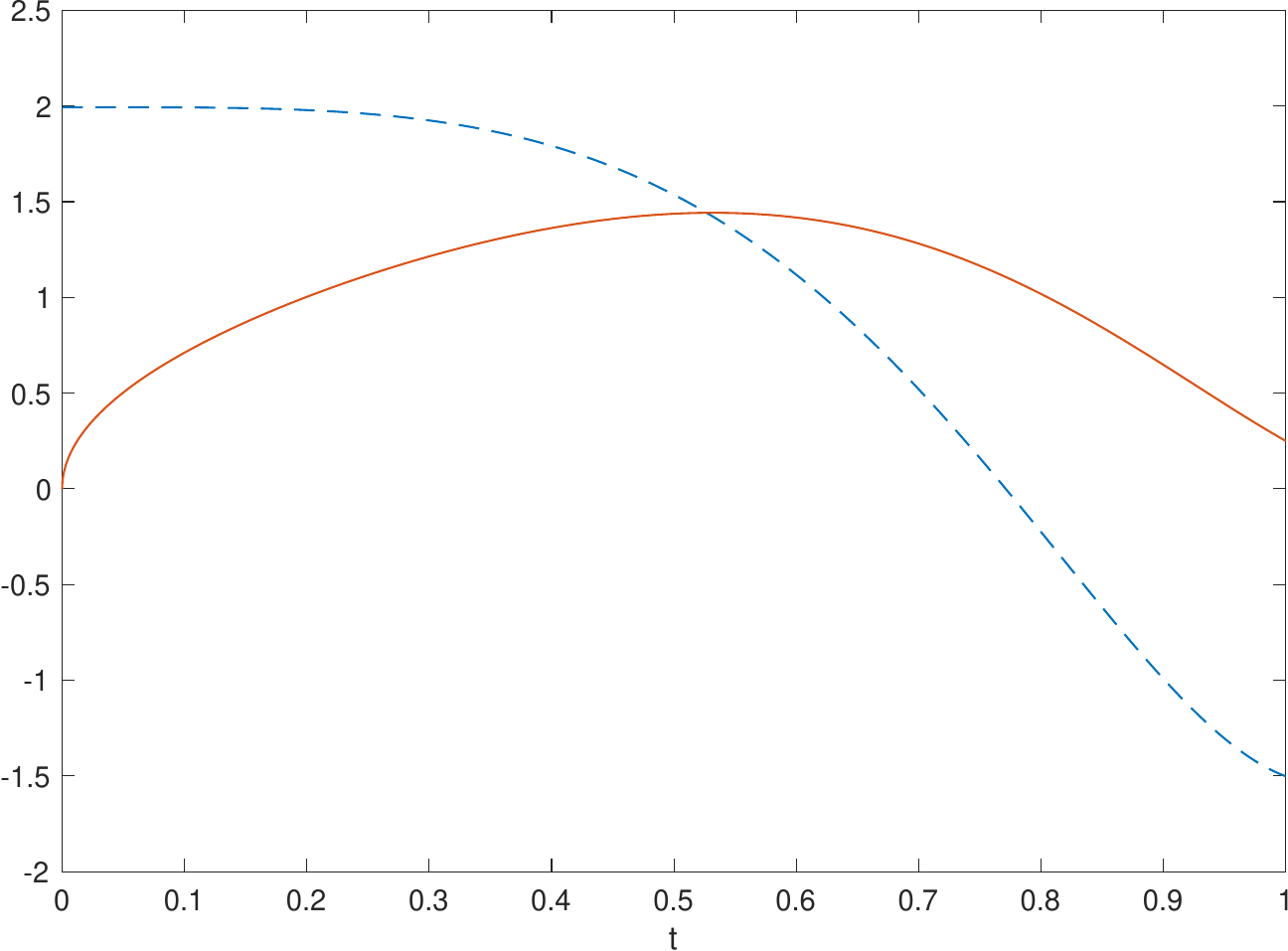}
\caption{solution (continuous line) and vector field (dashed line) for problem (\ref{prob2}).}
\label{prob1fig}
\end{figure}

The next problem that we consider is from \cite{DiFoFr2005} (see also \cite{Ga2018}):
\begin{equation}\label{prob2}
y^{(0.5)} = -y^{1.5} +\frac{40320}{\Gamma(8.5) }t^{7.5} -
           3\frac{\Gamma(5.25)}{\Gamma(4.75)}t^{3.75} +
           \left( \frac{3}2t^{0.25} - t^4 \right)^3 + \frac{9}4\Gamma(1.5), \quad t\in[0,1], \quad y(0)=0,
\end{equation}
whose solution is
\begin{equation}\label{prob2sol}
y(t) = t^8-3\,t^{4.25}+\frac{9}4\,t^{0.5}.
\end{equation}
According to \cite{Ga2018}, ``this problem is surely of interest because, unlike several other problems often proposed in the literature, it does not present an artificial smooth solution, which is indeed not realistic in most of the fractional-order applications.'' Despite this, a constant timestep $h=1/N$ turns out to be appropriate since, unlike the solution (\ref{prob2sol}), the vector field (\ref{prob2}) is very smooth at the origin, as one may see in Figure~\ref{prob1fig}. In Table~\ref{tab4} we list the maximum error by using different values of $N$: as is clear, by using $s\ge8$, only 32 steps are needed to gain full machine accuracy for the computed solution.

\begin{table}
\caption{\label{tab4} Maximum error for Problem (\ref{prob2}), constant timestep $h=1/N$.}

\smallskip
\centering

{\tt
\begin{tabular}{|r|r r r r r| }
\hline
 $s\,\backslash\,N$    & 2  & 4  & 8 & 16  & 32 \\
 \hline
 1~~~~ & 9.22e-01 & 5.65e-02 & 1.28e-02 & 1.35e-02 & 9.12e-03  \\ 
 2~~~~ & 7.48e-03 & 2.68e-03 & 5.15e-04 & 8.02e-05 & 1.91e-05  \\ 
 3~~~~ & 2.02e-03 & 1.96e-04 & 1.23e-05 & 2.04e-06 & 5.07e-07  \\ 
 4~~~~ & 2.29e-04 & 8.42e-06 & 2.72e-07 & 3.55e-08 & 3.70e-09  \\ 
 5~~~~ & 1.63e-05 & 3.52e-07 & 4.43e-09 & 3.44e-10 & 1.62e-11  \\ 
 6~~~~ & 7.61e-07 & 9.80e-09 & 6.57e-11 & 2.26e-12 & 1.47e-13  \\ 
 7~~~~ & 4.11e-08 & 3.71e-10 & 9.02e-12 & 3.46e-13 & 2.18e-14  \\ 
 8~~~~ & 1.24e-09 & 6.02e-11 & 1.87e-12 & 6.54e-14 & 4.22e-15  \\ 
 9~~~~ & 4.56e-10 & 1.44e-11 & 4.27e-13 & 1.65e-14 & 1.11e-15  \\ 
10~~~~ & 1.40e-10 & 4.40e-12 & 1.33e-13 & 4.77e-15 & 8.88e-16  \\ 
20~~~~ & 4.93e-14 & 1.33e-15 & 6.66e-16 & 6.66e-16 & 8.88e-16  \\ 
\hline
\end{tabular}}
\end{table}

Next, we consider the following problems, taken from \cite{Sa2023}:
\begin{equation}\label{prob3}
y^{(1/3)} = \frac{t}{10} \left[ y^3 - (t^{2/3}+1)^3\right] + \frac{\Gamma(5/3)}{\Gamma(4/3)} t^{1/3}, \qquad t\in[0,1], \qquad y(0) = 1.
\end{equation}
whose solution is $y(t) = t^{2/3}+1$, and
\begin{equation}\label{prob4}
y^{(1/3)} = \frac{1}3\left(y^3 - t^4\right) + \Gamma(7/3)t, \qquad t\in[0,1], \qquad y(0) = 0, 
\end{equation}
whose solution is $y(t) = t^{4/3}$.

We solve Problem (\ref{prob3}) by using $h_1=10^{-11}$ and $r=1.2$. In so doing 130 timesteps are needed to cover approximately the integration interval. The obtained results are listed in Table~\ref{tab5}. Also in this case, we obtain full accuracy starting from $s=8$.

\begin{table} 
\caption{\label{tab5} Maximum error for Problems (\ref{prob3}) and (\ref{prob34}), $r=1.2$ and $h_1=10^{-11}$.}

\medskip
\centering

{\tt
\begin{tabular}{|r|r|r|}
\hline
$s$ & error (\ref{prob3}) & error (\ref{prob34})\\
\hline
  1 & 3.25e-02  & ***\\ 
 2 & 8.86e-05  &5.13e-04\\ 
 3 & 8.36e-07  &4.21e-06\\ 
 4 & 1.41e-08  &7.55e-08\\ 
 5 & 3.03e-10  &1.63e-09\\ 
 6 & 7.54e-12  &3.95e-11\\ 
 7 & 3.46e-13  &1.06e-12\\ 
 8 & 2.09e-13  &2.09e-13\\ 
 9 & 2.09e-13  &2.09e-13\\ 
10 & 2.09e-13  &2.09e-13\\ 
20 & 2.09e-13  &2.09e-13\\ 

\hline
\end{tabular}}
\end{table}

In addition, we solve Problem (\ref{prob4}) by using a constant timestep $h=1/N$: in fact, the vector field can be seen to be a polynomial of degree 1. The obtained results are listed in Table~\ref{tab6}: an order 1 convergence can be observed for $s=1$ (which is consistent with (\ref{heq})), whereas full machine accuracy is obtained for $s\ge2$, due to the fact that, as anticipated above, the vector field of problem (\ref{prob4}) is a polynomial of degree 1 in $t$ and, consequently,  (\ref{inifrach}) and (\ref{sums0}) coincide, for all $s\ge2$.

\begin{table}
\caption{\label{tab6} Maximum error for Problem (\ref{prob4}), constant timestep $h=1/N$.}

\smallskip
\centering

{\tt
\begin{tabular}{|r|r r r r rr| }
\hline
 $s\,\backslash\,N$    & 2  & 4  & 8 & 16  & 32 & 64\\
 \hline
  1 &   *** & 1.56e-01 & 7.01e-02 & 3.59e-02 & 1.87e-02 & 9.75e-03 \\ 
 2 & 8.88e-16 & 1.33e-15 & 8.88e-16 & 8.88e-16 & 8.88e-16 & 8.88e-16 \\ 
 3 & 4.44e-16 & 6.66e-16 & 4.44e-16 & 4.44e-16 & 3.33e-16 & 4.44e-16 \\ 
 4 & 6.66e-16 & 6.66e-16 & 5.55e-16 & 2.22e-16 & 3.33e-16 & 5.55e-16 \\ 
 5 & 9.99e-16 & 9.99e-16 & 6.66e-16 & 5.55e-16 & 2.22e-16 & 5.55e-16 \\ 
 6 & 1.33e-15 & 7.77e-16 & 8.88e-16 & 6.66e-16 & 3.33e-16 & 5.55e-16 \\ 
 7 & 1.33e-15 & 8.88e-16 & 7.77e-16 & 4.44e-16 & 4.44e-16 & 7.77e-16 \\ 
 8 & 1.78e-15 & 1.11e-15 & 7.77e-16 & 6.66e-16 & 4.44e-16 & 6.66e-16 \\ 
 9 & 2.00e-15 & 1.22e-15 & 8.88e-16 & 8.88e-16 & 4.44e-16 & 6.66e-16 \\ 
10 & 2.00e-15 & 1.22e-15 & 8.88e-16 & 8.88e-16 & 4.44e-16 & 7.77e-16 \\ 
20 & 2.78e-15 & 1.89e-15 & 1.44e-15 & 1.11e-15 & 6.66e-16 & 8.88e-16 \\ 
\hline
\end{tabular}}
\end{table}

Finally, we reformulate the two problems (\ref{prob3}) and (\ref{prob4}) as a system of two equations, having the same solutions as above, as follows:
\begin{eqnarray}\nonumber
y_1^{(1/3)} &=& \frac{t}{10} \left[ y_1^3 - (y_2^{1/2}+1)^3\right] + \frac{\Gamma(5/3)}{\Gamma(4/3)} t^{1/3}, \qquad y_1(0) = 1,\\
\label{prob34}\\ \nonumber
y_2^{(1/3)} &=& \frac{1}3\left(y_2^3 - (y_1-1)^6\right) + \Gamma(7/3)t, ~~\qquad\qquad y_2(0) = 0, \qquad t\in[0,1]. 
\end{eqnarray}
We solve Problem (\ref{prob34}) by using the same parameters used for Problem (\ref{prob3}): $h_1=10^{-11}$, $r=1.2$, and 130 timesteps. The obtained results are again listed in Table~\ref{tab5}: it turns out that they are similar to those obtained for Problem (\ref{prob3}) and, also in this case, we obtain full accuracy starting from $s=8$.

\section{Conclusions}\label{fine}
In this paper we have devised a novel step-by-step procedure for solving fractional differential equations. The procedure, which generalizes that given in \cite{ABI2019}, relies on the expansion of the vector field along a suitable orthonormal basis, here chosen as the shifted and orthonormal Jacobi polynomial basis. The analysis of the method has been given, along with its implementation details. These latter details show that the method can be very efficiently implemented.  A few numerical tests confirm the theoretical findings.

It is worth mentioning that systems with FDEs of different orders can be also solved by slightly adapting the previous framework: as matter of fact, it suffices using different Jacobi polynomials, corresponding to the different orders of the FDEs. This, in turn, allows easily managing fractional differential equations of order $\alpha>1$, by casting them as a system of $\lfloor \alpha\rfloor$ ODEs, coupled with a fractional equation of order  $\beta\,:=\,\alpha-\lfloor\alpha\rfloor \in (0,1)$. 

As anticipated in Section~\ref{solve}, a future direction of investigation will concern the efficient numerical solution of the generated discrete problems, which is crucial when using larger timesteps. Also the optimal choice of the parameters for the methods will be further investigated,   as well as the possibility of adaptively defining the computational mesh. Further, the parallel implementation of the methods could be also considered, following an approach similar to \cite{AB1997}.



\end{document}